\documentclass[journal]{IEEEtran}
\usepackage{cite}

\usepackage[pdftex]{graphicx}

\usepackage{amsmath}
\usepackage{amsfonts}
\usepackage{amsthm}

\usepackage{amsmath}%
\usepackage{MnSymbol}%
\usepackage{wasysym}%

\theoremstyle{plain}
\newtheorem{thm}{Theorem}
\newtheorem{lem}{Lemma}

\newtheorem{remark}{Remark}

\newtheorem*{pf}{Proof}

\usepackage{algorithm}
\usepackage{algorithmic}
\usepackage{multirow}  
\usepackage{booktabs}  
\usepackage{amsthm}

\usepackage{array}

\ifCLASSOPTIONcompsoc
  \usepackage[caption=false,font=normalsize,labelfont=sf,textfont=sf]{subfig}
\else
  \usepackage[caption=false,font=footnotesize]{subfig}
\fi

\usepackage{stfloats}
\renewcommand{\arraystretch}{1.3}

\usepackage{nomencl}
\usepackage{ifthen}
\renewcommand{\nomgroup}[1]{%
    \ifthenelse{\equal{#1}{A}}{\item[\emph{\textbf{Set and Index}}]}{%
    \ifthenelse{\equal{#1}{B}}{\item[\emph{\textbf{Parameters}}]}{%
    \ifthenelse{\equal{#1}{C}}{\item[\emph{\textbf{Variables}}]}
    }
    }
    }
\makenomenclature

\begin{document}

\title{A Study on the Strong Duality of Conic Relaxation of AC Optimal Power Flow in Radial Networks}
%
%
%

\author{Xiaoyu Cao,~\IEEEmembership{Student Member,~IEEE,}
        Jianxue Wang,~\IEEEmembership{Senior Member,~IEEE,}
        and~Bo~Zeng,~\IEEEmembership{Member,~IEEE} \vspace{-25pt}}

\maketitle

\vspace{-3pt}\begin{abstract}
This letter proposes a set of closed-form conditions to ensure the strong duality of second-order cone program (SOCP) formulation for AC power flow in radial power networks. In addition, numerical evaluations on IEEE 33-bus test networks and a real-world distribution network are performed to demonstrate the validity of the proposed conditions.\vspace{-1pt}
\end{abstract} 

 \vspace{-10pt}

%
\IEEEpeerreviewmaketitle


\section{Introduction}
%
%
%
%
\vspace{-2pt}\IEEEPARstart{R}{elaxing} an  AC optimal power flow (OPF) model  into a convex program, e.g., second-order cone program (SOCP) or a semi-definite program (SDP)~\cite{jabr2006radial,bai2008semidefinite,lavaei2012zero,low2014convex}, enables us to utilize well-developed convex optimization tools to derive results stronger than those of the vastly adopted linear programming (LP) based DC power flow formulation. The SOCP formulation is particularly attractive as its computation burden is less heavy.
One great advantage of the LP formulation is its strong duality property, i.e., it  (the primal problem) shares the same optimal value as its dual linear program, which provides a substantial support to more advanced studies, e.g., developing  decomposition algorithms, deriving electricity market oriented decisions, and identifying the critical contingencies. Similarly, we note a few recent publications, e.g., \cite{D_Kekatos2014Stochastic,D_lee2015robust,D_lin2017decentralized,D_Haghighat2018bilevel,D_wu2018robust}, that critically rely on the strong duality of SOCP formulation.

However, it should be  pointed out that the strong duality of SOCP does not hold in general, i.e., a  non-zero duality gap may exist. Actually, we observe that the duality gap of some standard radial network could be non-negligible. Hence, theoretically speaking, if the strong duality cannot be guaranteed, results obtained based on the dual of SOCP formulation can only be considered as heuristic ones. In this letter,  we derive a set of closed form conditions on the network's physical parameters that guarantee the strong duality. \vspace{-5pt}

\section{Conic AC Optimal Power Flow Model}
\label{sec:OPF}
\begin{figure}[ht]
\centering
\includegraphics[width=2.4in]{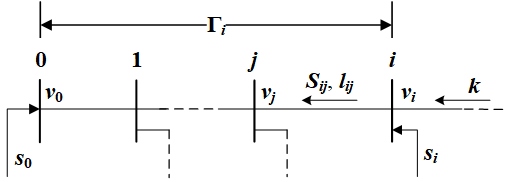} \vspace{-10pt}
\caption{Illustration of Radial Distribution Network} 
\label{fig:RN} \vspace{-5pt}
\end{figure}

Consider a distribution network with a radial topology, i.e., a spanning tree structure, as in Fig. \ref{fig:RN}. Let the radial network represented by $(\mathcal{N},\mathcal{E})$ where $\mathcal{E}$ denotes the set of branches and $\mathcal{N}=\{0\}\cup\{\mathcal{N}^+\}$ denotes the set of nodes that is the union of the substation node $0$ and the rest of nodes in  $\mathcal{N}^+=\{1,\ldots, n\}$. Given that the radial network has node $0$ as its root, node $i$ ($i\in \mathcal{N}^+$) has its unique parent node (denoted by $j$), and we can denote the unique path from node $i$ to node $0$ by $\Gamma_i$. Additionally, let $k$ to denote one child node of node $i$. \vspace{-5pt}

In the following, we present a branch flow model (BFM) based AC power flow formulation for this network. For this formulation, unless explicitly stated, the objective function $f$ is convex with no special structure. Note that the basic form of this formulation has been proposed in~\cite{Baran2002Optimal} and many variants have been adopted for different applications, e.g.,~\cite{low2014convex,D_Kekatos2014Stochastic,D_lin2017decentralized}.

\setlength{\arraycolsep}{0.0em}
\begin{eqnarray}
\vspace{-5pt}&&\mathbf{OPF}: \quad \min \ f({\rm Re}(\vec{s}),\vec{v},\vec{\ell})   \label{OBJ}\\
&&s.t. \ s_i=S_{ij}-\sum_{k:k\rightarrow i}(S_{ki}-z_{ki}\ell_{ki}), \quad \forall (i,j)\in \mathcal{E}\label{CONS1}\\
&&s_0=-\sum_{k:k\rightarrow 0}(S_{k0}-z_{k0}\ell_{k0}) \label{CONS2}\\
&&v_i-v_j=2{\rm Re}(\bar{z}_{ij} S_{ij})-|z_{ij}|^2\ell_{ij}, \quad  \forall (i,j)\in \mathcal{E} \label{CONS3}\\
&&\ell_{ij}=\frac{|S_{ij}|^2}{v_i}, \quad \forall (i,j)\in \mathcal{E}\label{CONS4}\\
&&0 \leq \ell_{ij} \leq \overline{\ell}_{ij}, \quad \forall (i,j)\in \mathcal{E}\label{CONS5}\\
&&\underline{v}_i \leq v_i \leq \overline{v}_i, \quad \forall i\in \mathcal{N}^+ \label{CONS6}\\
&&\underline{p}_i\leq {\rm Re}(s_i)\leq \overline{p}_i,\ \underline{q}_i\leq {\rm Im}(s_i)\leq \overline{q}_i, \quad \forall i\in \mathcal{N}^+ \label{CONS7}
\end{eqnarray}
Note that parameter $z_{ij}=r_{ij}+\mathbf{i}x_{ij}$ is the impedance of branch $(i,j)\in \mathcal{E}$ with $r_{ij}$ and $x_{ij}$ denoting the line resistance and reactance. Variable $\ell_{ij}$ denotes the squared magnitude of current on branch $(i,j)\in \mathcal{E}$, which is bounded by a positive number $\overline{\ell}_{ij}$ as in \eqref{CONS5}. Variable $v_i$ denotes the squared magnitude of voltage at node $i\in \mathcal{N}^+$, whose upper/lower bounds are denoted by $\overline{v}_i$/$\underline{v}_i$ as in \eqref{CONS6}. Variable $S_{ij}=P_{ij}+\mathbf{i}Q_{ij}$ is power flow though branch $(i,j)\in \mathcal{E}$, where $P_{ij}$ and $Q_{ij}$ denote the active and reactive power flows. Similarly, $s_i=p_i+\mathbf{i}q_i$ is the power injection at node $i\in \mathcal{N}^+$, where $p_i$ and $q_i$ denote the active and reactive power injections, and are bounded by $\overline{p}_i/\underline{p}_i$ and $\overline{q}_i/\underline{q}_i$ as in \eqref{CONS7}. The connection between nodal power injections and power flows is defined by branch flow equations \eqref{CONS1}-\eqref{CONS4}. The power injection at the substation node, i.e., node $0$, is $s_0=p_0+\mathbf{i}q_0$, and the associated square of the reference voltage level is fixed to constant $v_0$.

\vspace{-5pt}
\begin{remark}
\label{rmk:2}
Due to the nonlinear equality constraint \eqref{CONS4}, $\textbf{OPF}$ in \eqref{OBJ}-\eqref{CONS7} is non-convex. To convexify this formulation, \eqref{CONS4} is relaxed to the following second-order conic inequality~\cite{low2014convex}: 
\begin{equation}
\label{Constraint:SOCP1}\vspace{-2pt}
v_i\ell_{ij} \geq |S_{ij}|^2, \quad \forall (i,j)\in \mathcal{E}
\end{equation}
\noindent As a result, we obtain a convex relaxation of \textbf{OPF} defined by \eqref{OBJ}-\eqref{CONS3}, \eqref{CONS5}-\eqref{CONS7} and \eqref{Constraint:SOCP1}. We denote this relaxation by \textbf{OPF-Cr}. In particular, if the objective function $f$ is affine or convex quadratic,  \textbf{OPF-Cr} is an SOCP formulation (denoted by \textbf{OPF-SOCP}), which, under some sufficient conditions~\cite{low2014convex}, is exact and guarantees an optimal power flow solution to \textbf{OPF}.
\end{remark} \vspace{-5pt}

\vspace{-5pt}
\section{Sufficient Condition Ensuring Strong Duality}
Before theoretical derivations, we first make a few rather non-restrictive assumptions.
\begin{description}
  \item[A1.] The objective function $f$ is bounded from below.
  \item[A2.] The bounds on voltage magnitudes satisfy $\overline{v}_i>v_0>\underline{v}_i>0$ for all $i\in \mathcal{N}^+$. Typically, $\overline{v}_i$ and $\underline{v}_i$ are set within a small deviation around $v_0$. 
  \item[A3.]  The line resistance and reactance are positive, i.e., $r_{ij}>0$ and $x_{ij}>0$ for all $(i,j)\in \mathcal{E}$. \vspace{-10pt}
\end{description}

\subsection{ \textbf{OPF-Cr}  Reformulations with Restrictions}
To develop the sufficient conditions on strong duality, we construct two auxiliary SOCPs by reformulating constraints in \eqref{CONS1}-\eqref{CONS3} and \eqref{CONS5}-\eqref{Constraint:SOCP1}. First, by introducing new variables $\tau_{ij}\in \mathbb{R}_+$ and $\beta_{ij}\in \mathbb{R}$ for all $(i,j)\in \mathcal{E}$, we restrict our attention to a set of solutions of \textbf{OPF-Cr} (denoted by $\hat{S},\hat{\ell},\hat{s},\hat{v}$) that are represented as linear combinations of $\tau_{ij}$ and $\beta_{ij}$. Specifically,
\begin{eqnarray}
&&\hat{S}_{ij}=\frac{z_{ij}(\tau_{ij}-\beta_{ij})}{|z_{ij}|^2}, \quad \forall (i,j)\in \mathcal{E}  \label{PF:CONS1}\\
&&\hat{\ell}_{ij}=\frac{\tau_{ij}}{|z_{ij}|^2}, \quad \forall (i,j)\in \mathcal{E}  \label{PF:CONS2}
\end{eqnarray}
Plugging \eqref{PF:CONS1}-\eqref{PF:CONS2} into equalities \eqref{CONS1} and \eqref{CONS2}, variables $\hat{s}_i$ and $\hat{s_0}$ can be rewritten as
\setlength{\arraycolsep}{-0.2em}
\begin{eqnarray}
&&\hat{s}_i=\frac{z_{ij}}{|z_{ij}|^2}(\tau_{ij}-\beta_{ij})+\sum_{k:k\rightarrow i}\frac{z_{ki}}{|z_{ki}|^2}\beta_{ki}, \quad \forall i\in \mathcal{N}^+ \label{PF:CONS3}\\
&&\hat{s}_0=\sum_{k:k\rightarrow 0}\frac{z_{k0}}{|z_{k0}|^2}\beta_{k0}.  \label{PF:CONS4}
\end{eqnarray}
Similarly, equality \eqref{CONS3} can be converted into
\begin{equation}
\label{PF:CONS8A}
\hat{v}_i-\hat{v}_j=\tau_{ij}-2\beta_{ij}, \quad \forall (i,j)\in \mathcal{E}.
\end{equation}
Moreover, because the nodal voltage $\hat{v}_i$ in \eqref{PF:CONS8A} can be uniquely re-defined by summing right-hand-side (RHS) expressions over the connected path $\Gamma_i$, $\hat{v}_i$ can be rewritten as
\begin{equation}
\label{PF:CONS8}
\hat{v}_i=v_0+\sum_{(m,n)\in \Gamma_i} (\tau_{mn}-2\beta_{mn}), \quad \forall i\in \mathcal{N}^+.
\end{equation}

In addition to \eqref{CONS1}-\eqref{CONS3}, inequalities \eqref{CONS5}-\eqref{Constraint:SOCP1} should also be satisfied to ensure feasibility. Plugging \eqref{PF:CONS1}-\eqref{PF:CONS4} and \eqref{PF:CONS8} into conic constraint \eqref{Constraint:SOCP1} and affine constraints \eqref{CONS5}-\eqref{CONS7}, we have the first auxiliary SOCP (i.e., \textbf{OPF-SOCP$_1$}) as in \eqref{THM:CONS1}-\eqref{THM:CONS5}. 
\setlength{\arraycolsep}{-0.3em}
\begin{eqnarray}
&&(v_0+\!\sum_{(m,n)\in \Gamma_i}\!(\tau_{mn}-2\beta_{mn}))\tau_{ij} \geq \!(\tau_{ij}-\beta_{ij})^2, \forall (i,j)\in \mathcal{E} \label{THM:CONS1}\\
&&0 \leq \tau_{ij} \leq |z_{ij}|^2\overline{\ell}_{ij}, \quad \forall (i,j)\in \mathcal{E} \label{THM:CONS2}\\
&&\underline{v}_i-v_0 \leq \sum_{(m,n)\in \Gamma_i} (\tau_{mn}-2\beta_{mn}) \leq \overline{v}_i-v_0, \ \forall i\in \mathcal{N}^+  \label{THM:CONS3}\\
&&\underline{p}_i \leq \frac{r_{ij}}{|z_{ij}|^2}(\tau_{ij}-\beta_{ij})+\sum_{k:k\rightarrow i}\frac{r_{ki}}{|z_{ki}|^2}\beta_{ki} \leq \overline{p}_i, \ \forall i\in \mathcal{N}^+  \label{THM:CONS4}\\
&&\underline{q}_i \leq \frac{x_{ij}}{|z_{ij}|^2}(\tau_{ij}-\beta_{ij})+\sum_{k:k\rightarrow i}\frac{x_{ki}}{|z_{ki}|^2}\beta_{ki} \leq \overline{q}_i, \ \forall i\in \mathcal{N}^+  \label{THM:CONS5}
\end{eqnarray}

 Actually, by using the lower bound in \eqref{CONS6} as well as assumption A2, we can derive a further restriction on inequality \eqref{THM:CONS1} as in the following:
\begin{equation}
\label{THM:CONS6}
(\tau_{ij}-\beta_{ij})^2\leq \underline{v}_i\tau_{ij}, \quad \forall (i,j)\in \mathcal{E}
\end{equation}
With \eqref{THM:CONS6} being a conic inequality, our second auxiliary SOCP (i.e., {OPF-SOCP$_2$}) is defined by \eqref{THM:CONS2}-\eqref{THM:CONS5} and \eqref{THM:CONS6}.

\begin{remark}
\label{rmk:3}
As mentioned, auxiliary problems \textbf{OPF-SOCP$_1$} and \textbf{OPF-SOCP$_2$} define restricted solution spaces of \textbf{OPF-Cr}. We also note that \textbf{OPF-SOCP$_2$} is a restriction to \textbf{OPF-SOCP$_1$}. Hence, the feasible sets of \textbf{OPF-Cr}, \textbf{OPF-SOCP$_1$} and \textbf{OPF-SOCP$_2$} (denoted by $\mathbb{X}_{\rm OPF-Cr}$, $\mathbb{X}_{\rm OPF-SOCP_1}$ and $\mathbb{X}_{\rm OPF-SOCP_2}$) satisfy the following relationship: \vspace{-5pt}
\begin{equation}
\label{Relation}
\mathbb{X}_{\rm OPF-SOCP_2}\subset \mathbb{X}_{\rm OPF-SOCP_1} \subset \mathbb{X}_{\rm OPF-Cr}.\vspace{-15pt}
\end{equation}
\end{remark}

\subsection{Strong Duality of Reformulations and \textbf{OPF-Cr}} \vspace{-3pt}
Following \emph{Slater's Condition}, the strong duality holds for a general SOCP problem if either its primal problem or dual problem is bounded and strictly feasible, i.e., the problem is feasible and all the non-affine (conic) inequality constraints hold with strict inequalities \cite{ben2001lectures}. Hence, the \textbf{OPF-Cr} and our auxiliary SOCPs, which are bounded because of assumption A1, have the strong duality as long as  they are strictly feasible. Next, we consider the strict feasibility of \textbf{OPF-SOCP$_2$}.

\begin{lem}
\label{lem:1}
\textbf{OPF-SOCP$_2$} is strictly feasible and thus has the strong duality if any of the following conditions is satisfied:

\begin{description}
  \item[C1.] For every $i\in \mathcal{N}^+$, the bounds of its power injection satisfy either (i) $\underline{p}_i\leq0\leq \overline{p}_i, \ \underline{q}_i<0<\overline{q}_i$; or (ii) $\underline{p}_i<0<\overline{p}_i, \ \underline{q}_i\leq0\leq\overline{q}_i$; or (iii) $\underline{p}_i<0\leq\overline{p}_i, \ \underline{q}_i<0\leq\overline{q}_i$; or (iv) $\underline{p}_i\leq0<\overline{p}_i, \ \underline{q}_i\leq0<\overline{q}_i$.
  \item[C2.] $r_{ij}/x_{ij}\geq r_{ki}/x_{ki}$ for all $(i,j),(k,i)\in \mathcal{E}$; and $\underline{p}_i<0\leq\overline{p}_i, \ \underline{q}_i\leq 0\leq\overline{q}_i$ for all $i\in \mathcal{N}^+$.
  \item[C3.] $r_{ij}/x_{ij}\leq r_{ki}/x_{ki}$ for all $(i,j),(k,i)\in \mathcal{E}$; and $\underline{p}_i\leq 0\leq\overline{p}_i, \ \underline{q}_i<0\leq\overline{q}_i$ for all $i\in \mathcal{N}^+$.
\end{description}

\end{lem}

\begin{pf}
We make use of two new variables $\mu, \lambda_{ij} \in \mathbb{R}_{+}$ to simplify constraints over $\tau_{ij}$ and $\beta_{ij}$. Specifically, we have
\begin{eqnarray}
&&\tau_{ij}=\frac{|z_{ij}|^2 \overline{\ell}_{ij}}{\mu}, \quad \forall (i,j)\in \mathcal{E}  \label{PF:CONS10}\\
&&\beta_{ij}=\lambda_{ij} \tau_{ij}, \quad \forall (i,j)\in \mathcal{E}  \label{PF:CONS11}
\end{eqnarray}
Due to inequality \eqref{THM:CONS2}, we have $\mu \geq 1$. Then, plugging \eqref{PF:CONS10} and \eqref{PF:CONS11} into the strict version of inequality of \eqref{THM:CONS6} as well as affine inequalities \eqref{THM:CONS3}-\eqref{THM:CONS5}, we have
\begin{eqnarray}
&&\frac{(1-\lambda_{ij})^2}{\mu} < \frac{\underline{v}_i}{|z_{ij}|^2\overline{\ell}_{ij}}, \quad \forall (i,j)\in \mathcal{E} \label{PF:CONS12}\\
&&\underline{v}_i-v_0 \leq \sum_{(m,n)\in \Gamma_i} |z_{mn}|^2\overline{\ell}_{mn}\frac{1-2\lambda_{mn}}{\mu} \leq \overline{v}_i-v_0,\forall i\in \mathcal{N}^+ \label{PF:CONS13}\\
&&\underline{p}_i \leq r_{ij}\overline{\ell}_{ij} \frac{1-\lambda_{ij}}{\mu}+\frac{1}{\mu}\sum_{k:k\rightarrow i}r_{ki}\overline{\ell}_{ki}\lambda_{ki} \leq \overline{p}_i, \quad \forall i\in \mathcal{N}^+  \label{PF:CONS14}\\
&&\underline{q}_i \leq x_{ij}\overline{\ell}_{ij} \frac{1-\lambda_{ij}}{\mu}+\frac{1}{\mu}\sum_{k:k\rightarrow i}x_{ki}\overline{\ell}_{ki}\lambda_{ki} \leq \overline{q}_i, \quad \forall i\in \mathcal{N}^+  \vspace{-10pt}\label{PF:CONS15}
\end{eqnarray}
Clearly, with a fixed $\lambda_{ij}$ and $\mu\rightarrow +\infty$, the strict inequality in \eqref{PF:CONS12} can be easily achieved given that the left-hand-side of \eqref{PF:CONS12} approaches to $0+$ while its RHS is strictly positive following assumptions A2-A3. Also, through assumption A2, we have $\underline{v}_i-v_0<0<\overline{v}_i-v_0$. When $\mu\rightarrow +\infty$, inequality \eqref{PF:CONS13} holds as its middle term approaches to $0$.

Similarly, the feasibility of inequalities \eqref{PF:CONS14} and \eqref{PF:CONS15} can be ensured by properly setting the bounds of power injections, i.e., $\overline{p}_i/\underline{p}_i$ and $\overline{q}_i/\underline{q}_i$. let $\delta_{ij}^p,\delta_{ij}^q \in \mathbb{R}$ be such that
\begin{eqnarray}
&&\delta_{ij}^p = 1-\lambda_{ij}+\frac{\sum_{k:k\rightarrow i}\lambda_{ki}r_{ki}\overline{\ell}_{ki}}{r_{ij}\overline{\ell}_{ij}}, \quad \forall (i,j)\in \mathcal{E}  \label{PF:CONS16}\\
&&\delta_{ij}^q = 1-\lambda_{ij}+\frac{\sum_{k:k\rightarrow i}\lambda_{ki}x_{ki}\overline{\ell}_{ki}}{x_{ij}\overline{\ell}_{ij}}, \quad \forall (i,j)\in \mathcal{E}  \label{PF:CONS17}
\end{eqnarray}

\noindent Accordingly, inequalities \eqref{PF:CONS14}-\eqref{PF:CONS15} is equivalent to:
\begin{eqnarray}
&&\underline{p}_i \leq \frac{r_{ij}\overline{\ell}_{ij} \delta_{ij}^p}{\mu} \leq \overline{p}_i, \quad \forall i\in \mathcal{N}^+  \label{PF:CONS18}\\
&&\underline{q}_i \leq \frac{x_{ij}\overline{\ell}_{ij} \delta_{ij}^q}{\mu} \leq \overline{q}_i, \quad \forall i\in \mathcal{N}^+ \label{PF:CONS19}
\end{eqnarray}

\noindent Note that if $\delta_{ij}^p=0$ (or $\delta_{ij}^q=0$), the middle term of \eqref{PF:CONS18} or \eqref{PF:CONS19} will be zero so that it requires $\underline{p}_i\leq 0 \leq\overline{p}_i$ (or $\underline{q}_i\leq 0 \leq\overline{q}_i$) to ensure its feasibility; if $\delta_{ij}^p<0$ (or $\delta_{ij}^q<0$) and let $\mu\rightarrow +\infty$, the middle term of \eqref{PF:CONS18} or \eqref{PF:CONS19} will approach to $0-$, which requires $\underline{p}_i<0 \leq\overline{p}_i$ (or $\underline{q}_i<0 \leq\overline{q}_i$) to ensure its feasibility; if $\delta_{ij}^p>0$ (or $\delta_{ij}^q>0$), it requires $\underline{p}_i\leq 0<\overline{p}_i$ (or $\underline{q}_i\leq 0<\overline{q}_i$) to ensure its feasibility. For undetermined $\delta_{ij}^p$ and $\delta_{ij}^q$, constraints \eqref{PF:CONS18}-\eqref{PF:CONS19} are feasible as long as $\underline{p}_i<0<\overline{p}_i$ and $\underline{q}_i<0<\overline{q}_i$.

We also notice that for any $\lambda_{ki}$, either $\delta_{ij}^p$ or $\delta_{ij}^q$ can be specified by picking up a proper $\lambda_{ij}$, and then it dominates the value of the other.
So we can always find a $\lambda_{ij}$ to satisfy either (a) $\delta_{ij}^p=0$, or (b) $\delta_{ij}^q=0$, or (c) $\max\{\delta_{ij}^p, \delta_{ij}^q\}\leq 0$, or (d) $\min\{\delta_{ij}^p,\delta_{ij}^q\}\geq 0$ for every $(i,j)\in \mathcal{E}$, which corresponds to (i)-(iv) in condition \textbf{C1} that ensures the strict feasibility and thus the strong duality of \textbf{OPF-SOCP$_2$}.

Moreover, some alternative conditions of \textbf{C1} can be derived under special network parameters.
For instance, if $r_{ij}/x_{ij}\geq r_{ki}/x_{ki}$ for all $(i,j),(k,i)\in \mathcal{E}$, then we have
\setlength{\arraycolsep}{-0.4em}
\begin{eqnarray}
&&\frac{r_{ij}}{x_{ij}}\geq\frac{r_{ki}}{x_{ki}}\Leftrightarrow \frac{\sum_{k:k\rightarrow i}\lambda_{ki}r_{ki}\overline{\ell}_{ki}}{\sum_{k:k\rightarrow i}\lambda_{ki}x_{ki}\overline{\ell}_{ki}} \cdot \frac{x_{ij}\overline{\ell}_{ij}}{r_{ij}\overline{\ell}_{ij}}\leq 1 \nonumber\\
&&\Leftrightarrow \frac{\sum_{k:k\rightarrow i}\lambda_{ki}r_{ki}\overline{\ell}_{ki}}{r_{ij}\overline{\ell}_{ij}} \leq \frac{\sum_{k:k\rightarrow i}\lambda_{ki}x_{ki}\overline{\ell}_{ki}}{x_{ij}\overline{\ell}_{ij}} \Leftrightarrow \delta_{ij}^p \leq \delta_{ij}^q\label{Pf:C2}
\end{eqnarray}

\noindent Hence, there always exists a $\lambda_{ij}$ for every $(i,j)\in \mathcal{E}$ such that $\delta_{ij}^p\leq \delta_{ij}^q=0$. Then, it requires $\underline{p}_i<0\leq\overline{p}_i, \ \underline{q}_i\leq 0\leq\overline{q}_i$ as in \textbf{C2} to guarantee the strong duality of \textbf{OPFSOCP2}.

Similarly, when $r_{ij}/x_{ij}\leq r_{ki}/x_{ki}$, there exits a $\lambda_{ij}$ for every $(i,j)\in \mathcal{E}$ such that $\delta_{ij}^q\leq \delta_{ij}^p=0$, which requires $\underline{p}_i\leq 0\leq\overline{p}_i, \ \underline{q}_i< 0\leq\overline{q}_i$ as in \textbf{C3} to ensure the strong duality. $\hfill{} \blacksquare$
\end{pf}
Based on Lemma \ref{lem:1}, and the relationship among \textbf{OPF-SOCP$_2$}, \textbf{OPF-SOCP$_1$} and \textbf{OPF-Cr} (including the special case \textbf{OPF-SOCP}) stated in Remark \ref{rmk:3}, we have:
\begin{thm}
\label{thm:1}
\textbf{OPF-SOCP$_1$} and \textbf{OPF-Cr} are strictly feasible and thus have the strong duality when either of the conditions C1-C4 is satisfied.
\end{thm}

\section{Numerical Verification}
The proposed conditions are verified on IEEE 33-bus test networks and Southern California Edison (SCE) 56-bus network (a real-world distribution system). For illustration, some modifications are made on these networks to satisfy C1-C3 (e.g., by adding DGs or adjusting the network parameters). The primal/dual \textbf{OPF-SOCP}s of the original and modified networks are computed by MOSEK to attain the duality gaps. Table I shows the results of 1200 random instances of DG outputs. The average gap, the maximum gap, the number and ratio of instances with a gap less than $1e-4$ (which is a reasonable numerical standard to claim that the instance has the strong duality) are recorded in columns ``Avg-G", ``G$^+$", ``n\_SD", and ``r\_SD", respectively.

Note in Table I that: 1) Non-negligible duality gaps exist in many instances, e.g., 98.5\% of the instances in SCE 56-bus network with a maximum duality gap of 5.29\%. Such a large gap clearly cannot be attributed to numerical error, and therefore fails the strong duality test.  2) The strong duality of \textbf{OPF-SOCP} holds when some of conditions C1-C3 is satisfied. Our results support this claim as gaps of all 7,200 modified instances are numerically negligible.

Overall, we believe that the proposed closed-form conditions are significant that will support more advanced studies of SOCP formulation in radial networks. \vspace{-5pt}

\renewcommand\arraystretch{0.4}
\begin{table}[!h]
  \centering
  \caption{Tests on Different Distribution Networks} \vspace{-5pt}
  \setlength{\tabcolsep}{1.40mm}
  \scalebox{0.9}{
    \begin{tabular}{rlcccc}
    \toprule
    \multicolumn{2}{c}{Test Conditions} & Avg$-G$ & G$^+$ & N\_SD & R\_SD \\
    \midrule
          & Original System & 3.05E-03 & 6.99E-02 & 5     & 0.417\% \\
\cmidrule{2-6}    \multicolumn{1}{c}{IEEE 33-Bus } & Modified by C1 & 1.45E-08 & 1.54E-07 & 1200  & 100.0\% \\
\cmidrule{2-6}    \multicolumn{1}{c}{Network} & Modified by C2 & 2.43E-09 & 2.51E-08 & 1200  & 100.0\% \\
\cmidrule{2-6}          & Modified by C3 & 5.61E-09 & 5.91E-08 & 1200  & 100.0\% \\
    \midrule
          & Original System & 1.41E-02 & 5.29E-02 & 18     & 1.50\% \\
\cmidrule{2-6}    \multicolumn{1}{c}{SCE 56-Bus } & Modified by C1 & 8.00E-08 & 1.29E-05 & 1200  & 100.0\% \\
\cmidrule{2-6}    \multicolumn{1}{c}{Network} & Modified by C2 & 1.97E-06 & 4.52E-05 & 1200  & 100.0\% \\
\cmidrule{2-6}          & Modified by C3 & 1.94E-06 & 4.97E-05 & 1200  & 100.0\% \\
    \bottomrule
    \end{tabular}} \vspace{-15pt}
  \label{tab:sd}%
\end{table}

\ifCLASSOPTIONcaptionsoff
  \newpage
\fi



%




\bibliographystyle{IEEEtran}
\bibliography{test1}

\begin{thebibliography}{10}
\providecommand{\url}[1]{#1}
\csname url@samestyle\endcsname
\providecommand{\newblock}{\relax}
\providecommand{\bibinfo}[2]{#2}
\providecommand{\BIBentrySTDinterwordspacing}{\spaceskip=0pt\relax}
\providecommand{\BIBentryALTinterwordstretchfactor}{4}
\providecommand{\BIBentryALTinterwordspacing}{\spaceskip=\fontdimen2\font plus
\BIBentryALTinterwordstretchfactor\fontdimen3\font minus
  \fontdimen4\font\relax}
\providecommand{\BIBforeignlanguage}[2]{{%
\expandafter\ifx\csname l@#1\endcsname\relax
\typeout{** WARNING: IEEEtran.bst: No hyphenation pattern has been}%
\typeout{** loaded for the language `#1'. Using the pattern for}%
\typeout{** the default language instead.}%
\else
\language=\csname l@#1\endcsname
\fi
#2}}
\providecommand{\BIBdecl}{\relax}
\BIBdecl

\bibitem{jabr2006radial}
R.~A. Jabr, ``Radial distribution load flow using conic programming,''
  \emph{IEEE Trans. Power Syst.}, vol.~21, no.~3, pp. 1458--1459, 2006.

\bibitem{bai2008semidefinite}
X.~Bai, H.~Wei, K.~Fujisawa, and Y.~Wang, ``Semidefinite programming for
  optimal power flow problems,'' \emph{International Journal of Electrical
  Power \& Energy Systems}, vol.~30, no. 6-7, pp. 383--392, 2008.

\bibitem{lavaei2012zero}
J.~Lavaei and S.~Low, ``Zero duality gap in optimal power flow problem,''
  \emph{IEEE Trans. Power Syst.}, vol.~27, no.~1, pp. 92--107, 2012.

\bibitem{low2014convex}
S.~Low, ``Convex relaxation of optimal power flow-part \protect{I,II},''
  \emph{IEEE Trans. Control of Network Syst.}, vol.~1, 2014.

\bibitem{D_Kekatos2014Stochastic}
V.~Kekatos, G.~Wang, A.~J. Conejo, and G.~B. Giannakis, ``Stochastic reactive
  power management in microgrids with renewables,'' \emph{IEEE Trans. Power
  Syst.}, vol.~30, no.~6, pp. 3386--3395, 2014.

\bibitem{D_lee2015robust}
C.~Lee, C.~Liu, S.~Mehrotra, and Z.~Bie, ``Robust distribution network
  reconfiguration,'' \emph{IEEE Trans. Smart Grid}, vol.~6, pp. 836--842, 2015.

\bibitem{D_lin2017decentralized}
C.~Lin, W.~Wu, B.~Zhang, B.~Wang, W.~Zheng, and Z.~Li, ``Decentralized reactive
  power optimization method for transmission and distribution networks
  accommodating large-scale \protect{DG} integration,'' \emph{IEEE Trans.
  Sustain. Energy}, vol.~8, no.~1, pp. 363--373, 2017.

\bibitem{D_Haghighat2018bilevel}
H.~Haghighat and B.~Zeng, ``Bilevel conic transmission expansion planning,''
  \emph{IEEE Trans. Power Syst.}, 2018.

\bibitem{D_wu2018robust}
X.~Wu, A.~J. Conejo, and N.~Amjady, ``Robust security constrained
  \protect{ACOPF} via conic programming: \protect{Identifying} the worst
  contingencies,'' \emph{IEEE Trans. Power Syst.}, 2018.

\bibitem{Baran2002Optimal}
M.~Baran and F.~Wu, ``Optimal capacitor placement on radial distribution
  systems,'' \emph{IEEE Trans. Power Del.}, vol.~4, no.~1, pp. 725--734, 1989.

\bibitem{ben2001lectures}
A.~Ben-Tal and A.~Nemirovski, \emph{Lectures on modern convex optimization:
  \protect{Analysis}, algorithms, and engineering applications}.\hskip 1em plus
  0.5em minus 0.4em\relax SIAM, 2001, vol.~2.

\end{thebibliography}

%








\end{document}